\documentclass[12pt]{amsart}
\usepackage{amscd,amsmath,amssymb,amsfonts}
\theoremstyle{plain}
\newtheorem{thm}{Theorem}
\newtheorem{lem}[thm]{Lemma}
\newtheorem{cor}[thm]{Corollary}
\newtheorem{prop}[thm]{Proposition}

\theoremstyle{definition}
\newtheorem{defn}[thm]{Definition}

\newtheorem{rmk}[thm]{Remark}

\newtheorem{claim}[thm]{Claim}
\numberwithin{thm}{section}
\numberwithin{equation}{section}

\newcommand{\p}{\partial}

\newcommand{\ml}[2]{\begin{multline}\label{#1}#2 \end{multline}}
\newcommand{\ga}[2]{\begin{gather}\label{#1}#2 \end{gather}}

\newcommand{\surj}{\twoheadrightarrow}
\newcommand{\inj}{\hookrightarrow}

\newcommand{\Hom}{{\rm Hom}}

\newcommand{\Spec}{{\rm Spec \,}}

\newcommand{\sC}{{\mathcal C}}
\newcommand{\sD}{{\mathcal D}}
\newcommand{\sE}{{\mathcal E}}
\newcommand{\sF}{{\mathcal F}}

\newcommand{\sH}{{\mathcal H}}

\newcommand{\sL}{{\mathcal L}}

\newcommand{\sO}{{\mathcal O}}

\newcommand{\sQ}{{\mathcal Q}}

\newcommand{\sU}{{\mathcal U}}
\newcommand{\sV}{{\mathcal V}}
\newcommand{\sW}{{\mathcal W}}

\newcommand{\C}{{\mathbb C}}

\newcommand{\F}{{\mathbb F}}
\newcommand{\G}{{\mathbb G}}
\renewcommand{\H}{{\mathbb H}}

\renewcommand{\P}{{\mathbb P}}
\newcommand{\Q}{{\mathbb Q}}

\newcommand{\Rep}{{\rm Rep\hspace{0.1ex}}}
\newcommand{\Ind}{{\rm Ind\hspace{0.1ex}}}

\newcommand{\MIC}{\mbox{MIC}}
\newcommand{\id}{{\rm id}}
\newcommand{\stotimes}{\mbox{\hspace{.2ex}${\hspace{.2ex}}_s\otimes_t$\hspace{.2ex}}}
\newcommand{\sttimes}{\mbox{\hspace{.2ex}${\hspace{.2ex}}_s\times_t$\hspace{.2ex}}}
\begin{document}

\title[Gau{\ss}-Manin and Tannaka duality]{The Gau{\ss}-Manin connection and Tannaka duality}
\author{H\'el\`ene Esnault}
\address{Mathematik,
Universit\"at Duisburg-Essen, Mathematik, 45117 Essen, Germany}
\email{esnault@uni-essen.de}
\author{Ph\`ung  H\^o Hai}
\address{Mathematik,
Universit\"at Duisburg-Essen, Mathematik, 45117 Essen, Germany}
\email{ho-hai.phung@uni-essen.de}
\date{October 26, 2005}
\thanks{Partially supported by  the DFG Leibniz Preis.}
\begin{abstract}
If $f: X\to S$ is a fibration of complex connected analytic manifolds, given a point $x\in X$, there is an exact sequence of fundamental groups $0\to \pi_1(X_s, x)\to \pi_1(X,x) \to \pi_1(S, f(x))\to 0 $. Similarly, if $f: X \to {\rm Spec}(\F_q)$ is a smooth, geometrically connected  variety defined over a finite field, given a point $x\in X(\bar{\F}_q)$, there is an exact sequence of \'etale fundamental groups $0\to \pi_1(X\times_{\F_q} \bar{\F_q}, x)\to  \pi_1(X, x)\to
{\rm Gal}(\bar{\F}_q/\F_q)\to 0$. In particular, for any representation $\rho$ (resp. $\ell$-adic representation) of $\pi_1(X,x)$, one obtains an action of $\pi_1(S, f(x))$, resp.
${\rm Gal}(\bar{\F_q}/\F_q)$, on $H^i(\pi_1(X_s, x), \rho)$, resp.
$H^i(\pi_1(X\times_{\F_q} \bar{\F}_q, x), \rho)$. When  the natural map $H^i(\pi_1(X_s, x), \rho)\to H^i(X_s, \rho)$,
resp.
$H^i(\pi_1(X\times_{\F_q} \bar{\F}_q, x), \rho)\to
H^i(X\times_{\F_q} \bar{\F}_q, \rho)$ is an isomorphism,
this action defines the cohomology of the fibers as a local system over the base, resp. a Galois action on the cohomology of
$\rho$ over
$X\times_{\F_q} \bar{\F}_q$.

A good analog in algebraic geometry of the topological fundamental group on one side and
the \'etale fundamental group on the other side is provided by the Tannaka group associated to the category of flat connections. The action of the fundamental group or of the Galois group
of the base corresponds to the Gau{\ss}-Manin connection.
The purpose of this article is to show that the analogy isn't straightforward, that those actions are difficult to define, partly
because the homomorphism analogous to  $\pi_1(X_s, x)\to \pi_1(X,x)$ is not injective.

\end{abstract}
\maketitle
\begin{quote}

\end{quote}

\section{Introduction}
If $f: X \to {\rm Spec}(K)$ is a smooth, geometrically connected variety defined over a field of characteristic 0,
$K\supset k$ is a field extension, and $x\in X(K)$
is a rational point, one considers three Tannaka categories $\sC(X/K), \ \sC(X/k), \ \sC(K/k)$ of flat connections with compatible fiber functors. The objects of  $\sC(X/K)$ are  bundles (i.e. locally free coherent modules of finite type) with relative flat connections
$((\sV, \nabla_{/K}), \ \
\nabla_{/K}: \sV \to \Omega^1_{X/K}\otimes_{\sO_X} \sV)$, the ones of $\sC(X/k)$ are bundles with flat  absolute connections $((\sV, \nabla), \
\nabla: \sV \to \Omega^1_{X/k}\otimes_{\sO_X} \sV))$, the ones of $\sC(K/k)$ are
$K$-vector spaces with flat connections $((V,\nabla), \ \nabla: V\to \Omega^1_{K/k}\otimes_K V)$. The morphisms are the flat morphisms and the fiber functor has values in the category of finite
dimensional vector spaces ${\rm Vec}_K$ over $K$, defined by the restriction of $\sV$ to $x$ for $\sC(X/K), \ \sC(X/k)$ and by $V$ for $\sC(K/k)$.  Then $\sC(X/K)$ is a neutral Tannaka category, and Tannaka duality \cite{DeMil}, Theorem 2.11,  yields  the existence of a pro-group
scheme  $G(X/K)$ over $K$, so that $\sC(X/K)$ becomes equivalent to the representation category ${\rm Rep}(G(X/K))$. The two other Tannaka categories $\sC(X/k), \ \sC(K/k)$ are not
necessarily neutral. We assume that they are defined over $k$, which is to say that ${\rm End}_{\sC(K/k)}((K, d_{K/k}))=k.$
Then Tannaka duality \cite{DeGroth}, Th\'eor\`eme 1.12, yields the existence of  groupoid schemes
$G(X/k),  \ G(K/k)$ over $k$ acting on ${\rm Spec}(K)\times_k {\rm Spec}(K)$, so that, in the groupoid sense, $\sC(X/k)$ (resp. $\sC(K/k)$)
becomes equivalent to the representation category  ${\rm Rep}(K: G(X/k))$ (resp. $ {\rm Rep}(K: G(K/k)$). 
(We refer to the Appendix for the brief review of Deligne's theory of groupoid schemes.) 
Since the fiber functors are compatible, one has 
natural transformations $\sC(X/k)\xrightarrow{{\rm rest}} \sC(X/K)$ mapping an absolute connection $(\sV, \nabla)$ to the induced relative one $(\sV, \nabla/K)$ and $ \sC(K/k)\xrightarrow{f^*} \sC(X/k)$ 
mapping $(V,\nabla)$ to $f^*(V,\nabla)$. This yields
homomorphisms
\ga{1.1}{G(X/K)\xrightarrow{{\rm rest}} G(X/k)^\Delta,  \ G(X/k)\xrightarrow{f^*} G(K/k),}
where $^\Delta$ defines the induced group scheme over $K$ which is the restriction of 
$G(X/k)$ to the diagonal $\Delta={\rm Spec}(K)\to {\rm Spec}(K)\times_k {\rm Spec}(K)$, viewed as a group scheme over $K$.

On the other hand, if $V$ is an object in ${\rm Rep}(G(X/K))$, its cohomology group $H^i(G(X/K), V)$ is well defined \cite{Jan}, I, section 4, and is represented by an $i$-extension. Via the Tannaka formalism, this $i$-extension yields a $i$--extension of connections in \ $\sC(X/K)$ of the trivial connection $(\sO_X, d_{/K})$ by $(\sV, \nabla_{/K})$ corresponding to $V$. Via the connecting homomorphism $\delta: H^0_{DR}(X, (\sO_X, d_{/K}))\to H^i_{DR}(X, (\sV, \nabla_{/K}))$, where $_{DR}$ means de Rham cohomology relative to $K$, one defines a homomorphism
\ga{1.2}{H^i(G(X/K), V)\to H^i_{DR}(X, (\sV, \nabla_{/K})), \ i-{\rm extension}\mapsto \delta(1) .}
We show in Proposition \ref{prop2.2} that  \eqref{1.2} is an isomorphism for $i=0, 1$, is injective for $i=2$, thus in particular is an isomorphism when $X$ is an affine curve,
and also is an isomorphism if $X$ is a projective curve of genus $g\ge 1$.  

If $V$ is an object of   ${\rm Rep}(K: G(X/k))$, corresponding to the connection $(\sV, \nabla)$ and its restriction $(\sV, \nabla/K)$, then one has the Gau{\ss}-Manin connection defined on the finite dimensional $K$-vector space $H^i_{DR}(X, (\sV, \nabla/K))$. Via \eqref{1.2}, it corresponds to a groupoid action of $G(K/k)$ on $H^i(G(X/K), V)$ for $i=0,1$. We investigate the question of whether one can interpret the Gau{\ss}-Manin connection as one does in topology or in the $\ell$-adic categories as mentioned in the abstract above. In those categories, $G(X/K)$ corresponds to $\pi_1(X_s, s)$, resp. 
$\pi_1(X\times_{\F_q} \bar{\F}_q, x)$, $G(X/k)$ to $\pi_1(X,x)$, while $G(K/k)$ corresponds to $\pi_1(S, f(s))$, resp. ${\rm Gal}(\bar{\F}_q/\F_q)$. So via the standard exact sequence expressing the absolute fundamental group as an extension of the one on the base with the relative one, one defines an action of $\pi_1(S, f(x))$, resp. ${\rm Gal}(\bar{\F}_q/\F_q)$ on $H^i(\pi_1(X_s,x), \rho)$, resp. $H^i(\pi_1(X\times_{\F_q} \bar{\F}_q, x), \rho)$ for a representation $\rho$ of $\pi_1(X,x)$. 

In our setting, the map $f^*$ of groupoids in \eqref{1.1} is surjective, and one defines its kernel 
\ga{1.3}{0\to L^\delta\to G(X/k)\xrightarrow{f^*} G(K/k)\to 0}
as a discrete $k$-groupoid scheme, that is the $k$-linear morphism $L^\delta \to {\rm Spec}(K)\times_k {\rm Spec}(K)$ factors through the diagonal $\Delta \to {\rm Spec}(K)\times_k {\rm Spec}(K)$.
Its representation category  is equivalent to the representation category of the underlying group scheme $L\subset G(X/k)^\Delta$ over $K$. We show (Theorem \ref{thm5.8}) that  the representation category of $L$  is equivalent to the full  subcategory of $\sC(X/K)$, the objects of which are both subobjects and quotients of objects in $\sC(X/k)$. On the other hand, the homomorphism rest in \eqref{1.1} factors naturally through $L$. One defines the  subgroup scheme
\ga{1.4}{H:={\rm rest}(G(X/K))\subset L}
over $K$ of $L$ and show (Proposition \ref{prop3.1}) that its representation category is the full subcategory of $\sC(X/K)$ the objects of which are subquotients of objects in $\sC(X/k)$. This description allows us to define
an obstruction, local at $\infty$ of $X$, for an object of $\sC(X/K)$ to lie in ${\rm Rep}(H)$.  We show that this obstruction does  not necessarily vanish, thus
the homomorphism $G(X/K)\to G(X/k)^\Delta$ of group schemes is not necessarily injective (Proposition \ref{prop3.2}). This implies that
the kernel
of $G(X/K)\to G(X/k)^\Delta$ has a non-trivial subgroup scheme which is defined in categorial terms,  and which has the property that it has no nontrivial homomorphism into the additive group $\G_a$ (Theorem \ref{thm4.7}).

On the other hand, Deligne shows  that
 any relative subconnection of an absolute connection
is also the quotient as a relative connection of an absolute connection (Theorem \ref{thm5.10}). Thus the description of ${\rm Rep}(H) $ and ${\rm Rep}(L)$ as full subcategories of
${\rm Rep}(G(X/K))$ allows to conclude that $H=L$, that is the sequence
\ga{1.5}{G(X/K)\to G(X/k)\xrightarrow{f^*} G(K/k)\to 0}
is exact, where exacteness means that one sees the $k$-discrete groupoid scheme $L^\delta$ as a $K$-group scheme $L$, and then it is the image of $G(X/K)$ (Theorem \ref{thm5.11}).
\\ \ \\
Thus   there are at least two reasons why one can't overtake the standard argument describing a canonical action of $G(K/k)$ on $H^i(G(X/K), V)$, where $V$ is a finite representation of $G(X/k)$. Firstly, \eqref{1.1} is a sequence of groupoid schemes rather than group schemes, secondly $G(X/K)\to L$ is not injective. The reason why nevertheless one has this $G(K/k)$-groupoid action on $H^i(G(X/K), V)$ comes from the fact that
 the natural homomorphisms
$H^i(L, V)\to H^i(H, V)\to H^i(G(X/K), V)$ defined by functoriality are all isomorphisms for $i=0,1$ (Corollary \ref{cor4.3} and Theorem \ref{thm5.12}). 
As a corollary, one obtains a Tannaka theoritic formulation of the Gau{\ss}-Manin connection on $H^i_{DR}(X, (\sV, \nabla/K))$ (section 6). 

Finally, 
one can develop the same theory using connections relative to a finite 
dimensional $K$-vector space $D$ of $T_{K/k}$. More precisely,  let us remark that a connection is the same as an $\sO_X$-coherent module endowed with an action of the sheaf of differential operators. Let $D\subset T_{K/k}$ be a finite
dimensional  $K$-linear subspace of the tangent vectors of $K/k$, which is closed under brackets, such that ${\rm End}_D(K)=k$, and let $T_{X/k, D}$ be the inverse image of $f^*D$ in $T_{X/k}$, the  tangent sheaf of $K/k$, under the map $T_{X/k}\to f^*T_{K/k}$. Then $T_{X/k,D}$ is an extension of $f^*D$ by $T_{X/K}$, which generates a subalgebra $\sD_{X/k, D}\subset \sD_{X/k}$. The methods developed in this article could be used to treat the relation $\sC(X/K), \sC(X/k, D), \sC(K/k, D)$ as well, where
$\sC(X/k, D)$ is the category of $\sO_X$-coherent modules with an action of $\sD_{X/k, D}$,
and $\sC(K/k,D)$ is the category of finite dimensional $K$-vector spaces with a $\sD(K/k, D)$, the subalgebra of $\sD(K/k)$ spanned by $D$. We do not write the details.

\ \ \\

{\it Acknowledgements:} The first author started discussing in 2000 with Spencer Bloch on a possible Tannaka viewpoint on the Gau{\ss}-Manin connection. We thank him  for the fruitful exchanges we had at the time. We  thank  Nick Katz and Takeshi Saito for interesting discussions. We thank Pierre Deligne for his interest and for his  help. His comments allowed us to improve the results of an earlier version. 
Most specifically, Theorem \ref{thm5.10} is due to him.

\section{The neutral Tannaka category of flat connections}
Let $f:X \to {\rm Spec}(K)$ be a smooth geometrically connected variety defined over a field of characteristic 0.
\begin{defn} \label{defn2.1}
The category $\sC(X/K)$ of flat connections relative to $K$ (or simply of flat connections $/K$) has for objects the  flat connections $\big((\sV, \nabla), \nabla: \sV\to \Omega^1_{X/K}\otimes_{\sO_X} \sV\big)$, where $\sV$ is a locally free coherent module of finite type, and for morphisms 
the flat morphisms.
\end{defn}
It is a rigid abelian tensor category over $K$ and if we fix a $K$-rational point $x\in X(K)$, we can endow $\sC(X/K)$ with the fiber functor 
\ga{2.1}{\sC(X/K)\xrightarrow{\omega}  {\rm Vec}_K,\ (\sV, \nabla)\mapsto \omega((\sV, \nabla))=\sV|_x=:V}
with values in the category of finite dimensional $K$-vector spaces. Thus $\sC(X/K)$ becomes a neutral Tannaka category and by the fundamental Tannaka duality \cite{DeMil}, Theorem 2.11, 
$\omega$ defines an equivalence of tensor categories
\ga{2.2}{\sC(X/K)\xrightarrow{\omega \ \cong} {\rm Rep}_f(G(X/K))}
where $G(X/K)$ is the Tannaka group scheme over $K$, and ${\rm Rep}_f(G(X/K))$ is the category of its finite dimensional representations.  $G(X/K)$ is a pro-group scheme over $K$ which fulfills
\ml{2.3}{G(X/K)=\varprojlim_V G(V) \\ G(V)={\rm Im}(G(X/K)) \ {\rm in} \ GL(\omega((\sV, \nabla))).}
Let $V$ be an object of ${\rm Rep}_f(G(X/K))$. One defines its cohomology 
$H^i(G(X/K), V)$. 
Recall from \cite{Jan}, I, section 4 that if $G$ is a  group scheme, its cohomology is defined as the right derived functor of the functor $V\mapsto V^G$ of invariants, and is computed explicitely by cochains. Here $G(X/K)$ is pro-algebraic, acts on $V$ via its quotient $G(V)$, and the functor  
 $V\mapsto V^{G(X/K)}$ of invariants  factors through 
$V\mapsto V^{G(V)}=V^{G(X/K)}$. Setting $\sO[G(X/K)]:=\varinjlim_V \sO[G(V)]$ for the 
$K$-algebra of functions, with its canonical $G(X/K)$-action, the $G(X/K)$-injective modules are still direct summands of $({\rm trivial})\otimes_K \sO[G(X/K)]$ as in \cite{Jan}, I, (3.10). There are enough injective modules in this category of representations, and one defines  $H^i(G(X/K), V)$ as the right derived functor to the functor 
$V\mapsto V^{G(X/K)}$ of invariants. As $ V^{G(X/K)}={\rm Hom}_{G(X/K)}(K, V)$, 
one  has  as in \cite{Jan}, I, 4.2 (1) that cohomology is also the derived functor ${\rm Ext}^i_{G(X/K)}(K, V)$ to $V\mapsto 
{\rm Hom}_{G(X/K)}(K, V)$
\ga{2.4}{ H^i(G(X/K), V)={\rm Ext}^i_{{\rm Rep}_f(G(X/K))}(K, V).}
On the other hand, if $e$ is an $i$-extension of $K$ by $V$ in ${\rm Rep}_f(G(X/K))$, via Tannaka duality \eqref{2.2} one has an $i$-extension $\epsilon$ of $(\sO_X, d)$ by $(\sV, \nabla)$ in $\sC(X/K)$ with
\ga{2.5}{ \omega(\epsilon)=e} 
yielding a connecting homomorphism
\ga{2.6}{ \delta_{\epsilon}: H^0_{DR}(X, (\sO_X, d))\to H^i_{DR}(X, (\sV, \nabla)),}
where $ H^i_{DR}(X, (\sV, \nabla)):=\H^i(X, \Omega^\bullet_{X/K}\otimes_{\sO_X} \sV)$ is de Rham cohomology of the connection $(\sV, \nabla)$. 
This defines a  homomorphism of $K$-vector spaces
\ga{2.7}{\delta^i(X/K): H^i(G(X/K), V)\to H^i_{DR}(X, (\sV, \nabla)), \ e \mapsto \delta_{\epsilon}(1).}
\begin{prop} \label{prop2.2} 
The homomorphism $\delta^i(X/K)$  is an isomorphism for $i=0, 1$, and is injective for $i=2$. In particular, if $X/K$ is an affine curve, $H^2(G(X/K), V)=H^2_{DR}(X, (\sV, \nabla))=0$. Moreover, if $X/K$ is a smooth projective curve of genus $g\ge 1$, then it is an isomorphism for $i=2$.  
\end{prop}
\begin{proof}
For $i=0$, $H^0(G(X/K), V)\subset V$ is the largest trivial $G(X/K)$-subrepresentation. 
Thus by Tannaka duality \eqref{2.2}, it corresponds to $$(H^0_{DR}(X, (\sV, \nabla))\otimes_K \sO_X, 1\otimes d)$$ which is the largest trivial subconnection of $(\sV, \nabla)$, where $\omega((\sV, \nabla))=V$.  
For $i=1$, \eqref{2.4} says that a class $e\in H^1(G(X/K), V)$ is represented by an extension $e: 0\to V\to W\to K\to 0$ in ${\rm Rep}_f(G(X/K))$ and that two such extensions
$e, e'$ yield the same cohomology class if there is a commutative diagram
\ga{2.8}{\begin{CD}
0@>>> V @>>>  W @>>>  K @>>> 0\\
@. @V=VV @VVV @VV=V\\
0@>>> V @>>>  W' @>>>  K @>>> 0
\end{CD}
}
On the other hand, a class $\epsilon \in H^1_{DR}(X, (\sV, \nabla))$, with Cech cocycle $(u_{ij}, v_j)\in \sC^1(\sV)\times \sC^0(\Omega^1_{X/K}\otimes \sV), \delta(u)=\nabla(u)-\delta(v)=\nabla(v)=0$ on a Cech covering $\sU=\cup_i U_i$,  is represented by an extension $\epsilon: 
0\to (\sV, \nabla)\to (\sW, \nabla_{\sW}) \to (\sO_X,d)\to 0$ with $\sW|_{U_i}=(\sV\oplus \sO)|_{U_i},  \ \nabla_{\sW}(0\oplus 1)=v_i, \nabla_{\sW}|_{\sV \oplus 0}=\nabla$. Two such extensions $\epsilon, \epsilon'$ 
yield the same cohomology class if and only if there is a commutative diagram
\ga{2.9}{\begin{CD}
0@>>> (\sV,\nabla) @>>>  (\sW, \nabla_{\sW}) @>>>  (\sO_X,d) @>>> 0\\
@. @V=VV @VVV @VV=V\\
0@>>> (\sV,\nabla)  @>>> (\sW', \nabla_{\sW'})  @>>>  (\sO_X,d) @>>> 0
\end{CD}
}
Thus Tannaka duality \eqref{2.2} yields the result for $i=1$. 
We analyze now $i=2$. Let 
$e: 0\to V\to V_2\to V_1\to K\to 0$ be a 2-extension in ${\rm Rep}_f(G(X/K))$. We set $V_0={\rm Ker}(V_1\to K)$ and denote by $(\sV, \nabla)$ etc. the corresponding flat connections. To say that $\delta_{\epsilon}(1)=0$ is to say that there is a flat connection $(\sW, \nabla)$ 
which makes the following diagram a diagram with exact rows and columns
\ga{2.10}{
\begin{CD}
@. 0 @. 0\\
@. @AAA @AAA\\
0@>>> (\sV_0, \nabla) @>>> (\sV_1, \nabla) @>>> (\sO_X,d)@>>> 0\\
@. @AAA @AAA @AA=A\\
0 @>>> (\sV_2, \nabla) @>>> (\sW, \nabla) @>>> (\sO_X,d) @>>> 0\\
@. @AAA @AAA \\
@. (\sV, \nabla) @>>=> (\sV, \nabla) \\
@. @AAA @AAA \\
@. 0 @. 0
\end{CD}
}
But Tannaka duality yields the existence of an object $W$ in ${\rm Rep}_f(G(X/K))$ making the following diagram a diagram with exact rows and columns
\ga{2.11}{
\begin{CD}
@. 0 @. 0\\
@. @AAA @AAA\\
0@>>> V_0 @>>> V_1 @>>> K @>>> 0\\
@. @AAA @AAA @AA=A\\
0 @>>> V_2 @>>> W @>>> K @>>> 0\\
@. @AAA @AAA \\
@. V @>>=> V \\
@. @AAA @AAA \\
@. 0 @. 0
\end{CD}
}
This is to say that the class of $e$ in $H^2(G(X/K), V)$ dies. This shows injectivity.

We now prove the last part of the Proposition.
Let $X/K$ be a smooth projective curve, and $(\sV, \nabla)$ be a connection. Then $H^2_{DR}(X, (\sV, \nabla))$ is Poincar\'e dual to $H^0_{DR}(X, (\sV, \nabla)^\vee)$, where 
$(\sV, \nabla)^\vee$ is the dual connection. The inclusion 
$(H^0_{DR}(X, (\sV, \nabla)^\vee)\otimes_K \sO_X, 1\otimes d) \subset (\sV, \nabla)^\vee$ 
induces an isomorphism on $H^0_{DR}$, thus the dual projection
\ga{2.12}{(\sV, \nabla) \surj (H^2_{DR}(X, (\sV, \nabla))\otimes_K \sO_X, 1\otimes d)\cong 
\oplus_1^h (\sO_X, d)\\
 h={\rm dim}_K H^2_{DR}(X, (\sV, \nabla)) \notag.}
induces an isomorphism on $H^2_{DR}$. 
 On the other hand, assuming now that $g\ge 1$, 
there are two classes $\alpha, \beta \in H^1_{DR}(X)$, so that $0\neq \alpha\cup \beta \in H^2_{DR}(X)=K$.
Thus, there is a diagram of  extensions
\ga{2.13}{\begin{CD}
 (\sO_X, d) @>>> (\sE,\nabla) @>>> (\sO_X,d) \\
 @AAA \\
  (\sF,\nabla) \\
 @AAA \\
(\sO_X,d)
\end{CD}
}
in $\sC(X/K)$ which corresponds to $\alpha$ for the horizontal extension and $\beta$ for the vertical one. Denoting by  $(\sE,\nabla)_0$ the sub of $\oplus_1^h (\sE, \nabla)$ which is the inverse image of $(\sO_X,d)$ embedded diagonally in $\oplus_1^h (\sO_X,d)$, and setting $(\sF, \nabla)_0=\oplus_1^h (\sF, \nabla)$, \eqref{2.13} induces a 2-extension in $\sC(X/K)$
\ga{2.14}{0\to \oplus_1^h (\sO_X, d)\to (\sF, \nabla)_0\to (\sE, \nabla)_0 \to (\sO_X,d)\to 0}
which has the property that the image of connecting homomorphism $H^0_{DR}(X)\to \oplus_1^h H^2_{DR}(X)$ followed by a projection $\oplus_1^h H^2_{DR}(X)\to H^2_{DR}(X)$
is the fundamental class $\alpha\cup \beta$.
Since \eqref{2.12} induces an isomorphism on $H^2_{DR}$, there are connections
$(\sE_1,\nabla)$ and $(\sF_1, \nabla)$ together with a commutative diagram of 2-extensions
\ga{2.15}{\begin{CD}
0@>>> \oplus_1^h(\sO_X,d)@>>> (\sF, \nabla)_0 @>>> (\sE, \nabla)_0 @>>> (\sO_X, d) @>>> 0\\
@. @AAA @AAA @AAA @A=AA\\
0@>>> (\sV, \nabla)@>>> (\sF_1, \nabla) @>>> (\sE_1, \nabla) @>>> (\sO_X, d) @>>> 0
\end{CD}
}
in $\sC(X/K)$.
 We  apply Tannaka duality to the bottom 2-extension. This yields a 2-extension
\ga{2.16}{0\to V\to F_1\to E_1\to K\to 0}
in ${\rm Rep}_f(G(X/K))$, which has the property that the composite map
\ga{2.17}{H^0(G(X/K), K)\xrightarrow{{\rm connecting \ hom.}} H^2(G(X/K), V)
\xrightarrow{\delta^2(X/K)} \\
H^2_{DR}(X, (\sV, \nabla)) \xrightarrow{\eqref{2.12} \ {\rm iso}}
\oplus_1^h H^2_{DR}(X)\xrightarrow{{\rm proj.}} H^2_{DR}(X) \notag}
is an isomorphism.
This shows that $\delta^2(X/K)$ is surjective and finishes the proof.
\end{proof}

\begin{rmk}\label{rmk2.4}
If $K=\C$, then the classical Riemann-Hilbert correspondence establishes an equivalence of Tannaka categories
\ga{2.18}{{\rm Rep}_f(\pi_1^{{\rm top}}(X(\C), x))\cong \sC(X/\C).
}
It defines a homomorphism
\ga{2.19}{\pi_1^{{\rm top}}(X(\C), x) \to G(X/\C)(\C)}
with dense image. In particular, 
if $X/K=\P^1/\C$, then  $\pi_1^{{\rm top}}(X(\C), x)=0$, thus 
\ga{2.20}{H^2(G(X/K), K)=0 \ {\rm  while} \  H^2_{DR}(X,(\sO_X,d))=K.}
More generally, \eqref{2.17} gives a topological hint why the surjectivity on $H^2$ in Proposition \ref{prop2.2} is true  on a curve of genus $\ge 1$. Indeed in this case the universal covering of the underlying Riemann surface is contractible, thus the Hochshild-Serre spectral sequence degenerates and one has $H^2_{DR}(X/\C, (\sV, \nabla))=
H^2(\pi_1^{{\rm top}}(X(\C), x), V)$.

\end{rmk}
\section{The not necessarily neutral Tannaka category of flat connections}
Let $g: X\to {\rm Spec}(k)$ be a smooth  scheme with $k$ a field of characteristic 0, where $g$ factors through $f: X \to {\rm Spec}(K)$ as in Section 2, thus $X/K$ is a smooth
geometrically connected variety and  $K\supset k$ is a field extension with ${\rm End}_{\sC(K/k)}((K, d_{K/k}))=k$. We have $\sC(X/k)$ as in Definition \ref{defn2.1}.  The objects of $\sC(X/k)$ are flat connections $\big((\sV, \nabla), \ \nabla: \sV \to \Omega^1_{X/k}\otimes \sV\big)$, where $\sV$ is a coherent locally free sheaf on $X/K$, and morphisms are flat morphisms. This is a $k$-linear category. As in \eqref{2.1} we fix a $K$-rational point $x\in X(K)$ which defines a fiber functor
\ga{3.1}{\sC(X/k)\xrightarrow{\omega} {\rm Vec}_K, \ (\sV, \nabla)\mapsto \sV|_x=:V.}
Thus $\sC(X/k)$ becomes a non-neutral Tannaka category when $K\neq k$. By the fundamental
Tannaka duality \cite{DeGroth}, Th\'eor\`eme 1.12, there is a groupoid scheme $G(X/k)$ defined over $k$, acting on
${\rm Spec}(K)\times_k {\rm Spec}(K)$ so that $\omega$ defines an equivalence of tensor 
categories
\ga{3.2}{\sC(X/k)\xrightarrow{\omega \ \cong} {\rm Rep}_f(K:G(X/k))}
where $\Rep_f(K:G(X/k))$ denotes the category of finite dimensional $K$-representations
of $G(X/k)$. See the Appendix for a summary of the facts on groupoid schemes which will be used in the sequel.

We denote by $G(X/k)^\Delta$
the restriction of $G(X/k)$ to the diagonal.
Then $G(X/k)^\Delta$ is a discrete groupoid scheme over $k$.
 The representation category of a discrete groupoid scheme is equivalent to the representation category of the underlying
group scheme over $K$.
The embedding of Tannaka categories $\sC(X/k)\xrightarrow{{\rm rest}}\sC(X/K), \
(\sV, \nabla)\mapsto  {\rm rest}((\sV, \nabla))=(\sV, \nabla/K) $ with compatible fiber functor yields a homomorphism
\ga{3.3}{G(X/K)\xrightarrow{{\rm rest}} G(X/k)^\Delta \subset G(X/k).}
\begin{prop} \label{prop3.1}
The representation category on finite dimensional $K$-vector spaces of $H:={\rm rest}(G(X/K))\subset G(X/k)^\Delta$ is equivalent to the full subcategory of $\sC(X/K)$ the objects of which are subquotients of objects ${\rm rest}((\sV, \nabla))$.
\end{prop}
\begin{proof}
Let us denote by $\sC$ the  full subcategory of $\sC(X/K)$ the objects of which are subquotients of objects ${\rm rest}((\sV, \nabla))$, and by $G(\sC)$ its Tannaka group scheme.
Recall from \cite{DeMil}, Proposition 2.21 (a) that $G(X/K)\to G(\sC)$ is
faithfully flat  if and only if ${\rm rest}$ is fully faithful, which in our case is trivial, and
any subobject in $\sC(X/K)$ of an object in $\sC$ is an object in $\sC$, which is trivial as well in our case.
Recall from \cite{DeMil}, Proposition 2.21 (a) that $G(\sC)\to G(X/k)^\Delta$ is a closed immersion
if and ony if any object of $\sC$ is a subquotient of an object in ${\rm Rep}(G(X/k)^\Delta)$. But by definition, objects in $\sC$ are subquotients of objects in $\sC(X/k)$, thus as fortiori of objects
in ${\rm Rep}(\sC(X/k)^\Delta)$.
\end{proof}
\begin{prop} \label{prop3.2}
The homomorphism of group schemes  $G(X/K)\to G(X/k)^\Delta $ over $K$ is not necessarily injective.
\end{prop}
\begin{proof}
We assume that $X$ is an affine curve, and that $k=\C, K=\overline{\C(s)}$.
We wish to show that not every  connection on $X/K$ is a subquotient of a flat connection on $X/k$.
We consider a rank 1 connection $(\sL, \nabla)$ on $X/K$. Its formal completion
 at a point $y\in \bar{X}\setminus X$, which we assume to be $K$-rational with local parameter $t$,
is of the shape
\ml{3.4-8}{\sL\otimes_{\sO_X} K((t))=K((t))\cdot e,\\
 \ \big(\nabla\otimes_{\sO_X} K((t))\big)(e)= \alpha(t)\frac{dt}{t}\cdot e, \ \alpha(t)\in \frac{1}{t^n}K[[t]].}
Assume $(\sL,\nabla)$ is a subquotient of $(\sV, \nabla)$ on $X/k$.
By the Turritin-Levelt decomposition (see e.g. \cite{BBE}, section 5.9), one has
\ga{3.5}{(\sV,\nabla/K)\otimes_{\sO_X} K((t))=\oplus_i M_i\otimes U_i}
with $U_i$ nilpotent, $M_i$ irreducible, ${\rm Hom}(M_i, M_j)=K\cdot \delta_{ij}$.
Thus $(\sL,\nabla)\otimes_{\sO_X} K((t))$ has to be one of the $M_i$, say $M_0$, and not only it is a subquotient, it is also a sub relative connection.
We write the matrix of the connection in a basis adapted to the decomposition \eqref{3.5}
\ga{3.6}{\begin{pmatrix}
a &0\\
0 & b
\end{pmatrix}\frac{dt}{t} +
\begin{pmatrix}
A & B\\
C & D
\end{pmatrix}ds
}
with coefficients of $a,b, A,B,C,D$ in $ K((t))$, $a\frac{dt}{t}$ describing the connection $/K$ on $M_0\otimes U_0$, so $a=\alpha \otimes N_0$ with $N_0$ nilpotent with $K$-coefficients. The integrability condition implies
\ga{3.7}{\p_s a -t\p_t A=[a,A].}
This implies
\ga{3.8}{{\rm Tr \ Res}_{t=0} \p_s(a)\frac{dt}{t}=0}
as $[a, A]$ has trace zero and $(\p_t A) dt$ has  residue zero.
Let us denote by $\alpha_0 \in K$ the constant term in the $t$-expansion of $\alpha(t)\in \frac{1}{t^n}K[[t]]$. 
If
\ga{3.9}{\alpha_0\in K\setminus \C}
the condition \eqref{3.8} is not fulfilled and $(\sL, \nabla)\otimes_{\sO_X}K((t))$ is not a 
subconnection of
a flat connection $/\C$. Now starting with $(\sL,\nabla)$ on $X/K$, we can always achieve the condition
\eqref{3.9}. We possibly replace $X$ by a smaller affine $ X \setminus \Sigma$ so that
$\Gamma(\bar{X}\setminus \Sigma, \omega(y))$ contains a differential  form $\gamma$ so that $\alpha_0 +({\rm res}_{y} \eta) \in K\setminus \C$, and then replace $(\sL, \nabla)$ by
$((\sL, \nabla)\otimes (\sO_{X\setminus \Sigma}, d+\gamma))$
on $X\setminus \Sigma$.
 This finishes the proof.
\end{proof}
\begin{cor} \label{cor3.3} The kernel $N={\rm Ker}(G(X/K)\to G(X/k)^\Delta)$ is not trivial.
\end{cor}
We will show in Theorem \ref{thm4.7} that  $N$ has  a nontrivial subgroup with no $\G_a$ quotient.

\section{The universal de Rham extension}
In this section, the general assumption is as in Section 3: $g:X\to {\rm Spec}(k)$ is a smooth scheme
with $k$ a field of characteristic 0, with factorization $f: X\to {\rm Spec}(K)$ which makes $X$ a smooth geometrically connected
variety over the extension $K\supset k$. We assume that ${\rm End}_{\sC(K/k)}((K, d_{K/k}))=k$ and in this section, we assume in addition that
  the transcendence
degree of $K/k$ is $\le 1$. Fixing $x\in X(K)$, we have $\sC(X/k)$ and its groupoid scheme $G(X/k)$ as in \eqref{3.1}.

Let $(\sV, \nabla)$ be a flat connection on $X/k$. Recall \cite{KaNil}, section 3, that it Gau{\ss}-Manin connection is defined as the connecting homomorphism $H^i_{DR}(X, (\sV, \nabla/K))\xrightarrow{GM} \Omega^1_K\otimes_K H^i_{DR}(X, (\sV, \nabla/K))$
on 
the relative cohomology $H^i_{DR}(X, (\sV, \nabla/K)):=\H^i(X, \Omega^\bullet_{X/K}\otimes \sV)$ of the extension
\ga{4.1}{0\to \Omega^1_{K/k}\otimes_K (\Omega^{\bullet -1}_{X/K}\otimes_{\sO_X} \sV) \to
\Omega^\bullet_{X/k}\otimes_{\sO_X} \sV
\to \Omega^\bullet_{X/K}\otimes_{\sO_X} \sV\to 0. }
\begin{rmk}\label{rmk4.1}
As mentioned in \cite{KaNil}, (1.2), the de Rham cohomology $H^i_{DR}(X, (\sV, \nabla_{/K}))$ is also the right derived functor to the left exact global section functor 
\ga{4.2}{\MIC(X/K)\to {\rm Vec}_K^{qc}, \ (\sV, \nabla_{/K}) \mapsto H^0_{DR}(X, (\sV, \nabla_{/K})),} where $\MIC(X/K)$ is the category of quasi-coherent modules endowed with a flat connection on $X/K$, the morphisms being the flat morphisms, and ${\rm Vec}_K^{qc}$ is the category of quasi-coherent $K$-vector spaces (i.e. of infinite dimensional $K$-vector spaces). This category has enough injectives. 
One defines $\MIC(X/k)$ as the category of quasi-coherent modules endowed with a flat
connection on $X/k$. This category has enough injectives. 
Moreover, as tacitly mentioned in \cite{KaNil}, Remark 3.1,  
the restriction functor $\MIC(X/k)\to  \MIC(X/K), (\sV, \nabla)\mapsto (\sV, \nabla/K)$ sends  injectives to injectives. Indeed, the sheaf of algebras of differential operators $ 
{\rm PDDiff}(X/k)$ is flat over its sheaf of subalgebras $ {\rm PDDiff}(X/K)$
as Zariski locally one can choose coordinates. 
Therefore, 
the restriction functor has an exact left adjoint, that is for $M$ defined $/k$, with restriction to $K$ denoted by $M/K$, and $N$ defined  $/K$, one has
\ml{}{{\rm Hom}_{{\rm PDDiff}(X/k)}({\rm PDDiff}(X/k)\otimes_{{\rm PDDiff}(X/K)} N, M)
=\\
{\rm Hom}_{{\rm PDDiff}(X/K)}(N,M/K).\notag}
This implies that the restriction to $K$ of injective modules $/k$ are injective $/K$. 
 (We thank N. Katz for explaining us in more details his remark.) Thus 
 for $(\sV, \nabla)$ an absolute connection,  and $(\sV,\nabla)\xrightarrow{{\rm resolution}} I^\bullet$ an injective resolution in $\MIC(X/k)$, one has 
$$H^i_{DR}(X, (\sV, \nabla/K))=H^i( H^0_{DR}(X, I^\bullet/K)).$$ 
On the other hand, the restriction of $\nabla$ to 
$H^0_{DR}(X, (\sV, \nabla/K))\subset \sV$ is the Gau{\ss}-Manin connection which we denote by $GM$. This induces the Gau{\ss}-Manin connection on $H^i_{DR}(X, (\sV, \nabla/K))$, which we still denote by $GM$.
One obtains a commutative diagram
\ga{4.3}{\begin{CD}
\MIC(X/k)@> {\rm rest} >> \MIC(X/K)\\
@VH^i_{DR}(/K) VV @VV H^i_{DR} V\\
\MIC(K/k) @> {\rm rest} >> {\rm Vec}^{qc}_K
         \end{CD}
}
We still denote by $GM$ the Gau{\ss}-Manin connection on the full subcategory $\sC(X/k)\subset \MIC(X/k)$. Then \eqref{4.3} contains the sub-commutative square
\ga{4.4}{\begin{CD}
\sC(X/k)@> {\rm rest} >> \sC(X/K)\\
@VH^i_{DR}(/K) VV @VV H^i_{DR} V\\
\sC(K/k) @> {\rm rest} >> {\rm Vec}_K
         \end{CD}
}

\end{rmk}
\begin{thm} \label{thm4.2}
Let $(\sV, \nabla)$ be an object in $\sC(X/k)$. Then there is an extension in $\sC(X/k)$
\ga{4.5}{
0\to (\sV, \nabla)\to (\sW, \nabla)\to (H^1_{DR}(X, (\sV, \nabla/K))\otimes_K \sO_X, GM\otimes d)
\to 0
}
with the property that the connecting homomorphism 
\ml{4.6}{H^0_{DR}((H^1_{DR}(X, (\sV, \nabla/K))\otimes_K \sO_X, (GM\otimes d)/K))=H^1_{DR}(X, (\sV, \nabla/K))\\\xrightarrow{{\rm connecting}} H^1_{DR}(X, (\sV,\nabla/K))}
is the identity.
\end{thm}
\begin{proof}
The $K$-vector space $\sH om\big(H^1_{DR}(X, (\sV, \nabla/K)), H^1_{DR}(X, (\sV, \nabla/K))\big)$ is endowed with the connection $\sH om(GM, GM)$.
The identity $$1\in {\rm Hom}\big(H^1_{DR}(X, (\sV, \nabla/K)), H^1_{DR}(X, (\sV, \nabla/K))\big)$$
is a flat section. 
We apply \eqref{4.1} to the  flat connection 
$$ \sW:=\sH om(H^1_{DR}(X, (\sV, \nabla/K))\otimes \sO_X, \sV)$$
on $X/k$. Thus $1\in H^1_{DR}(X, \sW/K)$, with $GM(1)=0$. 
 Thus $1$ lifts to a class $$\tilde{1}
\in H^1_{DR}(X, \sW)=H^1_{DR}(X,  \sH om(H^1_{DR}(X, (\sV, \nabla/K))\otimes \sO_X, \sV)).$$
By the standard cocycle defined in the proof of Proposition \ref{prop2.2} for a class in 
$H^1_{DR}$, the class $\tilde{1}$ defines an extension \eqref{4.5} with \eqref{4.6} being the identity. 

\end{proof}
\begin{cor} \label{cor4.3}
Let $(\sV, \nabla)$ be an object in $\sC(X/k)$. Then the restriction homomorphism 
\ml{}{H^1(H, \omega((\sV,\nabla/K)))\to H^1(G(X/K), \omega(\sV, \nabla/K))=\\({\rm Proposition} \ \ref{prop2.2}) \ H^1_{DR}(X, (\sV, \nabla/K))\notag} is an isomorphism.

\end{cor}
\begin{proof}
As by Proposition \ref{prop3.1}, ${\rm Rep}_f(H)$ is equivalent to a full subcategory of $\sC(X/K)$, the homomorphism 
$$H^1(H, \omega((\sV,\nabla/K)))\to H^1(G(X/K), \omega(\sV, \nabla/K))$$ is injective. 
If now $e: 0\to (\sV, \nabla/K)\to (\sV', \nabla_{/K})\to (\sO_X,d/K)\to 0 $ is an extension in $\sC(X/K)$
with class $\bar{e}\in H^1_{DR}(X, (\sV, \nabla/K))$, then by Theorem \ref{thm4.2}, $e$ is isomorphic to the pull-back of \eqref{4.5} via $1\in \sO_X \mapsto \bar{e}\in H^1_{DR}(X, (\sV, \nabla/K))$, where \eqref{4.5} is now considered as an extension of relative connections on $X/K$. Consequently, $(\sV', \nabla_{/K})$ is a subconnection of an absolute flat connection, thus 
$(\sV', \nabla_{/K})$ is an object of ${\rm Rep}_f(H)$. This shows surjectivity and finishes the proof. 

\end{proof}
\begin{defn}\label{defn4.4} We define ${\rm Rep}_f(H)\subset {\rm Rep}_f(H)^t \subset {\rm Rep}_f(G(X/K))$ as the smallest full subcategory of ${\rm Rep}_fG(X/K)$ containing ${\rm Rep}_f(H)$ and which is thick, that is so that whenever two objects are in ${\rm Rep}_f(H)^t $, so is any extension. 
\end{defn}
\begin{lem} \label{lem4.5}
${\rm Rep}_f(H)^t$ is defined by its objects. One successively defines ${\rm Obj}_n$ as extensions in $
{\rm Rep}_f(G(X/K))$  of objects in ${\rm Obj}_{n-1}$, with ${\rm Obj}_0={\rm Obj}({\rm Rep}_f(H)))$.  Then the objects of ${\rm Rep}_f(H)^t$ consist of the union of the ${\rm Obj}_n$. 
\end{lem}
\begin{proof}
We just have to see that the category we construct in this manner is a sub-Tannaka category. It is obviously stable by taking duals. Furthermore, if $V$ is in ${\rm Obj}_n$ and $W\subset V$, then 
$W\in  {\rm Obj}_n$. Indeed, write $0\to V_1\to V\to V_0\to 0$ a defining extension for $V$ with $V_i\in {\rm Obj}_{n-1}$, then we just just have to see that $V_1\cap W$ is still in ${\rm Obj}_{n-1}$. Thus by induction on $n$, we just have to see that if $V\in {\rm Obj}({\rm Rep}_f(H))$, and $W\subset V$ is a subobject in ${\rm Rep}_f(G(X/K)) $, then it is a subobject in ${\rm Rep}_f(H)$. This is to say that $G(X/K)\surj H$, which is the definition. Finally, the definition is trivially compatible with tensor products. 
\begin{rmk} \label{rmk4.6} The example of Proposition \ref{prop3.2} has rank 1 thus shows that 
if $H^t$ is the Tannaka group of ${\rm Rep}_f(H)^t$, then one has
\ga{4.7}{G(X/K)\to H^t \to H}
with $G(X/K)\to H^t$ not necessarily injective. However it is  fully faithful as, as in the proof of Proposition \ref{prop3.1}, every subobject in ${\rm Rep}_f(G(X/K))$ of an object in ${\rm Rep}_f(H)^t$ is in
${\rm Rep}_f(H)^t$.
\end{rmk}
We define the non-trivial algebraic $K$-group
\ml{4.8}{N^t:={\rm Ker}(G(X/K)\surj H^t)\subset N=
{\rm Ker}(G(X/K)\surj H)\\\subset N={\rm Ker}(G(X/K)\surj H)
.}
\end{proof}
One has
\begin{thm} \label{thm4.7}
One has $H^1(N^t,\G_a)=0$.
\end{thm}
\begin{proof}
We consider the exact sequence 
\ga{4.9}{0\to N^t\to G(X/K)\to H^t\to 0}
of $K$-algebraic groups.  We consider the relative connection $(\sO_X, d/K)$. As it is the resriction of an absolute connection, $\omega((\sO_X, d))$ is certainly a representation of $H$.  From \eqref{4.7}, one has a factorization 
\ga{4.10}{
H^1(H, \omega((\sO_X,d)))\to H^1(H^t, \omega((\sO_X,d))) 
\to H^1(G(X/K), \omega((\sO_X,d))).}
Thus  Corollary \ref{cor4.3} allows to conclude that $$H^1(H^t, \omega((\sO_X,d))) \to H^1(G(X/K), \omega((\sO_X,d)))$$ is surjective. In fact it is injective as well as ${\rm Rep}_f(H^t)\subset {\rm Rep}_f (G(X/K))$ is a full subcategory. We don't use the injectivity. We conclude from the long exact sequence associated to \eqref{4.9} that $H^1(N^t, K)={\rm Hom}(N^t, \G_a)$ is equal to the kernel
of $H^2(H^t, K)\to H^2(G(X/K), K)$. ($K$ is here the trivial object). But this kernel is 0 again as 
${\rm Rep}_f(H^t)\subset {\rm Rep}_f(G(X/K))$ is full.
\end{proof}
\begin{rmk} \label{rmk4.8}
In order to get rid of the assumption
on the transcendence degree of $K/k$ being $\le 1$, one has to introduce the category of vertical connections, that is those connections on $X/k$, the curvature of which lies in $\Omega^2_K\otimes {\rm End}(\sV)$. This is because the universal extension \eqref{4.5} is a priori only a vertical connection, which is not necessarily flat. This introduces correspondingly the Tannaka groups $H_v, H_v^t$ etc. with the same conclusions as in Corollary \ref{cor4.3} and Theorem \ref{thm4.7}. Since we do not see any applications of this, we do not detail the construction.
\end{rmk}
\begin{rmk} \label{rmk4.9} If $X/K$ is an affine curve, then the embedding ${\rm Rep}_f(H)\xrightarrow{\iota} {\rm Rep}_f(G(X/K))$ is  thick. It is equivalent to saying that $H^i(H, V)\xrightarrow{\iota^*} H^i(G(X/K), V)$ is an isomorphism for $i=1$ and for all objects in ${\rm Rep}_f(H)$, which we show now.  An object of ${\rm Rep}_f(H)$
is of the shape $V=V'/V''$ with $V''\subset V'\subset W$ with $W=\omega((\sW, \nabla))$ and $(\sW, \nabla)$ an object in $\sC(X/k)$. By fullness (Proposition \ref{prop3.1}), $H^i(H, V)\xrightarrow{\iota^*} H^i(G(X/K), V)$ is injective. Thus applied to $V''$ and $i=2$ for which $H^2(G(X/K), V'')=0$ by Proposition \ref{prop2.2}, and to $i=1$ for $V'$ and $V$, we see that it is enough to show that $\iota^*$ is an isomorphism for $i=1$ and $V'$. Set $V_0=W/V' \in {\rm Obj}{\rm Rep}_f(H))$. Then again $\iota^*$ is injective for $V_0$ and $i=1$ while it is an isomorphism for $W/K$ and $i=1$ by Corollary \ref{cor4.3}. Using again injectivity if $\iota^*$ for $V'$ and $i=2$, we obtain 
the result.
\end{rmk}
\section{The exact sequence of groupoids}
The assumptions in this section are the same as in Section 3: $g:X\to {\rm Spec}(k)$ is a smooth scheme
with $k$ a field of characteristic 0, with factorization $f: X\to {\rm Spec}(K)$ which makes $X$ a smooth geometrically connected
variety over the extension $K\supset k$. We assume
${\rm End}_{\sC(K/k)} ((K, d_{K/k}))=k$.
Fixing $x\in X(K)$, we have $\sC(X/k)$ and its groupoid scheme $G(X/k)$ as in
Section 3. See also the Appendix.

We recall that $G(X/k)^\Delta$ is the discrete groupoid scheme, pull back of $G(X/k)$ over the diagonal $\Delta\to {\rm Spec}(K)\times_k {\rm Spec}(K)$. We define similarly $G(K/k)^\Delta$, pull back of $G(K/k)$ over the diagonal
$\Delta\to {\rm Spec}(K)\times_k {\rm Spec}(K)$. Both $G(X/k)^\Delta$ and $G(K/k)^\Delta$, as $K$-schemes,  are $K$-algebraic groups. 
\begin{lem} \label{lem5.1}
The   homomorphism of groupoid schemes  $$G(X/k)\xrightarrow{f^*} G(K/k)$$ is surjective, 
and induces a surjective homomorphism $$G(X/k)^\Delta \xrightarrow{f^*} G(K/k)^\Delta$$ of algebraic groups.
\end{lem}
\begin{proof}
The composite map $G(K/k)\xrightarrow{x^*} G(X/k)\xrightarrow{f^*} G(K/k)$ is the identity.
\end{proof}
We define the  pro-group scheme over $K$
\ga{5.1}{L={\rm Ker}(G(X/k)^\Delta\xrightarrow{f^*} G(K/k)^\Delta).}
Since the composite of functors
\ga{5.2}{\sC(K/k)\xrightarrow{f^*} \sC(X/k) \xrightarrow{{\rm rest}} \sC(X/K)}
sends any object to a finite sum of the trivial objet, the composite map of groupoid schemes
\ga{5.3}{G(X/K)\xrightarrow{{\rm rest}} G(X/k)\xrightarrow{f^*} G(K/k)}
fulfills
\ga{5.4}{{{\rm rest}}(G(X/K))=H\subset L.}
\begin{lem}[Key Lemma] \label{lem:key}
Let $V \in {\rm Obj}({\rm Rep}_f(K:G(X/k)))$.
Then $$H^0(L, V)= H^0(G(X/K), V).$$
\end{lem}
\begin{proof}
It is clear that 
\ga{5.5}{H^0(L, V)\subset H^0(G(X/K), V).}
 We wish to show surjectivity. 
We first show
\begin{claim} \label{claim5.3}
$$f^*H^0_{DR}(X, (\sV, \nabla/K))=
H^0_{DR}(X, (\sV, \nabla/K))
\otimes_K \sO_X \subset \sV $$
is the largest subbundle which is stabilized by $\nabla$ and on which  
\ml{5.6}{\nabla=f^*\delta \ {\rm with } \ 
\delta=\nabla|_{H^0_{DR}(X, (\sV, \nabla/K))} \  {\rm so}\\ 
(H^0_{DR}(X, (\sV, \nabla/K)), \delta)\in {\rm Obj}(\sC(K/k)).}
\end{claim}
\begin{proof}
Indeed, by flatness of $\nabla$, the composite map
\ga{5.7}{H^0_{DR}(X, (\sV, \nabla/K))\xrightarrow{\nabla} \Omega^1_K\otimes \sV\xrightarrow{\nabla} \Omega^1_K\otimes \Omega^1_{X/K}\otimes \sV}
is vanishing. On the other hand, one has
\ga{5.8}{\Omega^1_K\otimes\sV\xrightarrow{\nabla = 1_{\Omega^1_K}\otimes (\nabla/K)}
\Omega^1_K\otimes \Omega^1_{X/K}\otimes 
 \sV.}
One concludes 
\ml{5.9}{\nabla(H^0(X, (\sV, \nabla/K)))\subset \\
\Omega^1_K\otimes_K {\rm Ker}(\nabla/K)=\Omega^1_K\otimes H^0(X, (\sV, \nabla/K)).
}
Consequently, $H^0(X, (\sV, \nabla/K))\otimes_K \sO_X \subset \sV$
is stabilized by $\nabla$ and lies in the largest subbundle on which $\nabla$ is of the shape $f^*\delta$. On the other hand, it has to be the largest such, as any other $\sW\subset \sV$ would have the property that $(\nabla/K)|_{\sW}$ is generated by flat sections.

\end{proof}
Claim \ref{claim5.3} shows that $H^0(G(X/K), V)$ is a $G(X/k)$-representation on which $G(X/k)$ acts via its quotient $G(K/k)$, thus $H^0(G(X/K), V)\subset H^0(L, V)$. 
This finishs the proof. 
\end{proof}
\begin{cor} \label{cor5.4}
If $V_i, i=1,2$ are objects of ${\rm Rep}_f(K: G(X/k))$, with restrictions $V_i/K$ as objects of ${\rm Rep}_f(G(X/K))$,  
then $${\rm Hom}_L(V_1, V_2)={\rm Hom}_{G(X/K)}(V_1, V_2).$$
\end{cor}
\begin{proof}
We just set $V=V_1^\vee \otimes V_2$ and apply ${\rm Hom}(V_1, V_2)=H^0(V)$ in the corresponding category together with Lemma \ref{lem:key}.

\end{proof}
\begin{lem} \label{lem5.5}
Any finite dimensional representation of $G(X/k)^\Delta $ can be
embedded into the image of an object in $\sC(X/k)$ under
the fiber functor and consequently can be 
represented as a quotient of the image of an object in
$\sC(X/k)$ under the fiber functor.

\end{lem}
\begin{proof}
Any representation of $G(X/k)^\Delta $ can be embedded
in a direct sum of copies of the function algebra $\sO(G(X/k)^\Delta)$.
Indeed, the action of $G(X/k)^\Delta$ on a vector space
$V$ induces a coaction of the Hopf algebra $\sO(G(X/k)^\Delta)$
on $V$:
\ga{5.10}{\delta: V\to V\otimes_K \sO(G(X/k)^\Delta).}
The axiom for $\delta$ is precisely  requiring that
$V$ can be embedded into $\sO(G(X/k)^\Delta )^{\dim_KV}$
as a $G(X/k)^\Delta$-representation.

By construction, the function algebra is the image under the
fiber functor of an object in the Ind-category of the category
$\sC(X/k)$ (which consists of filtered direct limits of objects of
$\sC(X/k)$).  Thus one can find an object of $\sC(X/k)$, the
image of which contains the image of $V$ in
$\sO(G(X/k)^\Delta)^{\dim_KV}$.

The second claim is obtained by taking dual
representations.
\end{proof}
\begin{lem}\label{lem5.6}
Any representation of $L$ can be embedded into the restriction
to $L$ of representation of $G(X/k)^\Delta $ and consequently
can also be represented as the quotient of the restriction of a  representation
of $G(X/k)^\Delta $.
\end{lem}
\begin{proof}
Let 
\ga{5.11}{{\rm Ind} : {\rm Rep}(L) \to {\rm Rep}(G(X/k)^\Delta)}
be the right adjoint functor to the restriction functor
\ga{5.12}{{\rm Res}:  {\rm Rep}(G(X/k)^\Delta)\to {\rm Rep}(L).}
One has the  functorial isomorphism
\ga{5.13}{{\rm Hom}_L({\rm Res}( V),W)\xrightarrow{\cong} 
\Hom_{G(X/k)^\Delta}(V,{\rm Ind}( W)).}
It is well-known that $\sO(G(X/k)^\Delta)$ is faithfully flat over
it subalgebra $\sO(G(K/k)^\Delta)$ (\cite{waterhause}, Chapter 13).
This implies that the functor ${\rm Ind}$ is faithfully exact
(\cite{takeuchi79}, Chapter 2).
Setting $V={\rm Ind}( W)$ in \ref{5.13}, one obtains a canonical map
$u_W:{\rm Ind}( W)\to W$ which is non-zero whenever $W$ is non-zero,
since ${\rm Ind}$ is faithfully exact. We want to show that this
map is always surjective.

Let $U={\rm Im}(u_W)$ and $T=W/U$. We have the following diagram
\ga{5.14}{\begin{CD}
0&@>>>& {\rm Ind}( U) &@>>> &{\rm Ind}( W) &@>>>&{\rm Ind}( T) &@>>>&0\\
&&&&@VVV&&@VVV&&@VVV\\
0&@>>>&U&@>>>&W&@>>>&T&@>>>&0
\end{CD}}
The composition ${\rm Ind}( W)\surj {\rm Ind}( T)\to T$ is 0, therefore
${\rm Ind}( T)\to T$ is a zero map, implying $T=0$.

Since ${\rm Ind} (W)$ is a union of its finite dimensional subrepresentations,
we can therefore find a finite dimensional 
$G(X/k)^\Delta$-subrepresentation of 
${\rm Ind} (W)$ which still maps surjectively on $W$.

\end{proof}
As any infinite dimensional representation is the inductive limit of finite dimensional ones, one obtains 
\begin{cor} \label{cor5.7}
Every finite dimensional 
representation of $L$ can be embedded into the restriction to
$L$ of a finite dimensional representation of $G(X/k)$, consequently, it can
be represented as a quotient of the restriction to $L$
of a finite dimensional representation of $G(X/k)$.
\end{cor}
We are now in the position to prove the 
\begin{thm} \label{thm5.8}
The category ${\rm Rep}_f(L)$ is equivalent
to the full subcategory $\sC$ of $\sC(X/K)$, the objects of which are subobjects
as well as quotient objects of the restriction to $K$ of an absolute connection.  

\end{thm}
\begin{rmk} \label{rmk5.9}
The objects of $\sC$ are of the shape ${\rm Im}(\varphi)$ where $\varphi \in {\rm Hom}(V_1/K, V_2/K)$, with $V_i=\omega((\sV_i,\nabla_i))$,  $(\sV_i, \nabla_i)$
objects of $\sC(X/k)$ and $\varphi$ a morphism in $\sC(X/K)$.
\end{rmk}
\begin{proof}[Proof of Theorem \ref{thm5.8}]
Let us first remark that by definition, $\sC$ is a full subcategory
of $\sC(X/K)$, which is trivially closed under taking the tensor product.
We don't know yet whether it is an abelian subcategory.

We denote by  $\sQ: {\rm Rep}_f(L)\to \sC(X/K)$ the functor defined by
the homomorphism $q: G(X/K)\to L$. By
Corollary \ref{cor5.7} the image of $\sQ$ lies in $\sC$.

 Being a tensor functor, $\sQ$ is faithful. We show
that it is also full.
Let $U_0, U_1$ be objects  in Rep$_f(L)$ and $\phi:\sQ(U_0)\to \sQ(U_1)$
 a $\sC$-morphism,
i.e. $\phi$ is a $K$-linear map $U_0\to U_1$,
 which is only $G(X/K)$-linear where the
 actions of $G(X/K)$ is induced from the homomorphism $q:G(X/K)\to L$.
It is to show that $\phi$ is in fact $L$-linear.

 By Corollary \ref{cor5.7} there are
 $L$-linear morphisms $\pi:V_0\surj U_0$ and 
 $\iota:U_1\hookrightarrow V_1$, where $V_0$ and $V_1$ are objects in ${\rm Rep}_f(K:G(X/k))$.
One has
$U_i=\omega((\sU_i, \nabla_i))$ for relative connections $(\sU_i, \nabla_i) \in {\rm Obj}(\sC(X/K))$, $V_i=\omega((\sV_i, \nabla_i))$ for absolute connections
$(\sV_i, \nabla_i)\in {\rm Obj}(\sC(X/k))$
with
$(\sV_{0},\nabla_0/K)\xrightarrow{\pi \ \surj} (\sU_0,\nabla),
(\sU_1,\nabla_1)\xrightarrow{\iota \ \inj} (\sV_1,\nabla_1/K)$.
We set $\psi=\iota\phi\pi$.
\ga{5.15}{\begin{CD}U_0& @>\phi>> & V_0\\
@A\pi \ {\rm surj} AA&& @VV\iota \ {\rm inj}  V\\
V_0& @>>\psi> & V_1
\end{CD}}
By Corollary \ref{cor5.4}, the map $\psi$ is $L$-linear. This implies that $\phi$ is $L$-linear as well.

We now show that each
object of $\sC$ is isomorphic to the image under $\sQ$ of 
a representation of $L$. An object of $\sC$ has the form ${\rm Im}(\varphi)$
where $\varphi: V_0/K\to V_1/K$ as in Remark \ref{rmk5.9}.
By 
the above discussion,  $\varphi$ is also in the image of
$\sQ$, hence so is ${\rm Im}(\varphi)$.

Thus the functor $\sQ:{\rm Rep}_f(L)\to \sC$  is fully faithful
and each object of $\sC$ is isomorphic to the image of
an object of ${\rm Rep}_f(L)$. This shows that $\sC$ is a tensor subcategory
in $\sC(X/K)$, hence a Tannaka category, and the induced map
$\bar{q}: G(\sC)\to L$ is an isomorphism (\cite{DeMil}, Theorem 2.21).
\end{proof}
On the other hand, one has the following
\begin{thm}[Deligne] \label{thm5.10}
Let $(\sL, \nabla)$ be an object in $\sC(X/K)$ and assume there is an object $(\sV, \nabla_{\sV})\in \sC(X/k)$ so that in $\sC(X/K)$, one has an injection $(\sL, \nabla)\subset (\sV, \nabla_{\sV}/K)$. Then there is
an object $(\sW, \nabla_{\sW}) \in \sC(X/k)$ so that in $\sC(X/K)$, one has a surjection
$(\sW, \nabla_{\sW}/K)\surj (\sL, \nabla)$.
\end{thm}
\begin{proof}(Variant of Deligne's proof) We first assume that $(\sL, \nabla)$ is of rank 1. Then we define the $(\sL,\nabla)$-isotypical component $(\sW, \nabla)$ of $(\sV, \nabla_{\sV}/K)$ as follows. Set $(\sV',\nabla)=(\sV, \nabla_{\sV}/K)\otimes (\sL, \nabla)^\vee$, which is an object in $\sC(X/K)$. Then the inclusion $(\sL, \nabla)\subset (\sV, \nabla_{\sV}/K)$ corresponds to a nontrivial section in $H^0_{DR}(X, (\sV', \nabla))$. Set $\sV_1=\sV'/\big(H^0_{DR}(X, (\sV',\nabla))\otimes \sO_X\big)$ with induced connection
relative to $K$. If $H^0_{DR}(X, (\sV_1,{\rm induced \ connection/K}))=0$ then one defines $(\sW, \nabla)= H^0_{DR}(X, (\sV', \nabla))\otimes (\sL, \nabla)$. If not, define $\sV'_1$ to be the inverse image of $H^0_{DR}(X, (\sV_1, \nabla)) \otimes \sO_X$ in $\sV'$. The relative connection  on $\sV'$ induced from $\nabla_{\sV}/K$  stabilizes $\sV'_1$. Define $\sV_2=\sV'/\sV'_1$ with its induced connection. If $H^0_{DR}(X, (\sV_2,\nabla))=0$ one defines
similarly as before   $(\sW,\nabla)= (\sV'_1, {\rm induced \ connection})\otimes (\sL, \nabla)$. If not, we go on. So $(\sW, \nabla)$ is the largest subconnection of $(\sV, \nabla/K)$ which is a successive extension of $(\sL, \nabla)$ by itself. In particular, $(\sL, \nabla)$ is a quotient of $(\sW, \nabla)$ as well.

On the other hand, as $\nabla_{\sV}/K$ stabilizes $\sW$, one has
 $\nabla_{\sV}(\sW)\subset (\Omega^1_{X/k}\otimes_{\sO_X} \sV)' \subset \Omega^1_{X/k}\otimes_{\sO_X} \sV$, where $(\Omega^1_{X/k}\otimes_{\sO_X} \sV)' $ 
denotes the inverse image of 
$\Omega^1_{X/K}\otimes_{\sO_X} \sW$ via the projection $\Omega^1_{X/k}\otimes_{\sO_X} \sV \to \Omega^1_{X/K}\otimes_{\sO_X} \sV$. Consequently, the composite map $\sW\to \Omega^1_{X/k}\otimes_{\sO_X}\sV\to  \Omega^1_{X/k}\otimes_{\sO_X}(\sV/\sW)$, which is $\sO_X$-linear, has values in $f^{-1}(\Omega^1_{K/k})\otimes_{K}(\sV/\sW)$. We denote by 
$\ell: \sW\to f^{-1}(\Omega^1_{K/k})\otimes_{K} (\sV/\sW)$ the composite map. Integrability of $\nabla_{\sV}$ implies that the diagram
\ga{5.16}{\begin{CD}
\sW @>\ell>> f^{-1}(\Omega^1_{K/k})\otimes_{f^{-1}(K)} (\sV/\sW)\\
@V\nabla_{\sV}/K VV @VV 1_{\Omega^1_{K/k}}\otimes (\nabla_{\sV}/K )|_{(\sV/\sW)} V\\
\Omega^1_{X/K}\otimes_{\sO_X}\sW@>\ell>> f^{-1}(\Omega^1_{K/k})\otimes_{f^{-1}(K)} \Omega^1_{X/K}\otimes_{\sO_X} (\sV/\sW)
          \end{CD}
}
is commutative.

The right vertical map is the tensor product of the $K$-vector space $\Omega^1_{K/k}$ with the 
relative connection $(\sV/\sW) \xrightarrow{(\nabla_{\sV}/K )|_{(\sV/\sW)}} 
\Omega^1_{X/K}\otimes_{\sO_X}(\sV/\sW)$. Thus 
$\big(f^{-1}(\Omega^1_{K/k})\otimes_{f^{-1}(K)} (\sV/\sW), 1_{\Omega^1_{K/k}}\otimes (\nabla_{\sV}/K )|_{(\sV/\sW)} \big)$ is an object in the Ind-category spanned by $\sC(X/K)$, or simply  in $\sC(X/K)$ if   the tanscendence degree of $K/k$ is finite, and $\ell$ is a morphism in this category. Consequently $\ell(\sW)$ is a subobject 
in $\sC(X/K)$. This implies that $\ell(\sL)$ is a subobject in $\sC(X/K)$ as well, and the projection ${\rm pr}\circ \ell(\sL)$  in all  $\big((\sV/\sW), (\nabla_{\sV}/K )|_{(\sV/\sW)}\big)$ 
obtained by $K$-linear projection $\Omega^1_{K/k}\to K$ is a subobject in $\sC(X/K)$ as well. 
By definition of $\sW$, this implies that ${\rm pr}\circ \ell(\sL)=0$ for all such projections, thus $\ell(\sL)=0$, thus $\ell(\sW)=0$. 
 We conclude that $\nabla_{\sV}$ stabilizes $\sW$. We define $\nabla_{\sW}=\nabla_{\sV}|_{\sW}$. This finishes the proof in this case.

It remains to consider the case when $(\sL,\nabla)$ has higher rank $r$. We write
\ga{5.17}{ (\sL, \nabla)= {\rm det}(\sL, \nabla)\otimes \wedge^{r-1}(\sL, \nabla)^\vee,}
which shows the existence of a surjective map
\ga{5.18}{(\sW', \nabla'/K)\otimes\wedge^{r-1}(\sV, \nabla_{\sV}/K)^\vee \surj (\sL,\nabla)\subset (\sV, \nabla_{\sV}/K)}
where $(\sW', \nabla')$ is the object in $\sC(X/k)$  constructed in the rank 1 case for ${\rm det}(\sL, \nabla)$. We set $(\sW, \nabla_{\sW})= (\sW', \nabla')\otimes\wedge^{r-1}(\sV, \nabla_{\sV})^\vee$. 
This finishes the proof.

\end{proof}
Theorems \ref{thm5.8} and \ref{thm5.10},
together with  Proposition \ref{prop3.1}, allow us now to
conclude
\begin{thm} \label{thm5.11}
Assume as usual that
${\rm End}_{\sC(K/k)}((K, d_{K/k}))=k$. Then $H=L$, that is the sequence of groupoid schemes
\ga{5.19}{G(X/K)\to G(X/k)\xrightarrow{f^*} G(K/k)\to 0}
is exact, in the sense that one sees the $k$-groupoid scheme $L^\delta={\rm Ker}(f^*)$ as a $K$-group scheme and then ${\rm Im}(G(X/K))= {\rm Ker}(f^*).$
\end{thm}
\begin{proof}
Since both categories ${\rm Rep}_f(H) \supset {\rm Rep}_f(L)$ are  full
subcategories of ${\rm Rep}_f(G(X/K))$, it is enough to identify their objects.
If $V$ is an object of ${\rm Rep}_f(H)$, then there are $E_1\subset E\subset F$ in ${\rm Rep}_f(G(X/K)$, with $V=E/E_1$ and $F=\omega(\sF, \nabla/K)$, with $(\sF, \nabla)$ an object in $\sC(X/k)$.  By Theorem \ref{thm5.10}, the relative connection $(\sE, \nabla)$ with
$\omega((\sE, \nabla))=E$ is the quotient in $\sC(X/K)$ of an obsolute connection, so
the connection $(\sV, \nabla)$ with $V=\omega((\sV, \nabla))$ is the quotient in $\sC(X/K)$ of an obsolute connection as well. Applying this result to $V^\vee$, one concludes that $(\sV, \nabla)$ is also a subconnection in $\sC(X/K)$ of an object in $\sC(X/k)$. This finishes the proof.

\end{proof}

We now prove the theorem  which was one motivation for the article.
\begin{thm} \label{thm5.12}
Assume as usual that ${\rm End}_{\sC(K/k)}((K, d_{K/k}))=k$ and that the transcendence degree of $K/k$ is $\le 1$. 
Let $V\in {\rm Obj}({\rm Rep}_f(K:G(X/k)))$.
Then
$$H^i(L, V)=H^i(G(X/K), V) \ {\rm for} \ i=0, 1.$$
\end{thm}
\begin{proof}
This is an immediate consequence of  Theorem \ref{thm5.11} together with Corollary \ref{cor4.3}.

\end{proof}
\begin{cor} \label{cor5.13} Under the assumptions of Theorem \ref{thm5.12},
let $V=\omega((\sV, \nabla))$ be an object of ${\rm Rep}_f(K: G(X/k))$. Then
the extension defined in \eqref{4.5} yields an extension in ${\rm Rep}_f(K: G(X/k))$
\ga{5.20}{0\to V \to \omega(\sW, \nabla)\to H^1_{DR}(X, (\sV, \nabla/K))=H^1(L, V)\to 0}
with the property that the connecting homomorphism
\ga{5.21}{H^0(L, H^1(L, V))=H^1(L, V)\xrightarrow{{\rm connecting}} H^1(L, V)}
is the identity. 
\end{cor}

\section{The Gau\ss-Manin connection from the Tannaka viewpoint}
As usual, we consider an absolute connection $(\sV, \nabla)\in {\rm Obj}(\sC(X/k))$ together with its fiber functor $V=\sV|_x\in {\rm Obj}({\rm Vec}_K)$. The finite dimensional $K$-vector space
$H^0(L, V)$ is a $G(K/k)$-representation in a natural way. Indeed, for $(a,b): T\to {\rm Spec}(K)\times_k {\rm Spec}(K)$, and $g_{ab}\in G(K/k)(T)_{ab}$, consider $\tilde{g}_{ab}\in G(X/k)(T)$ a preimage. Then (see  Appendix, A.3) $\tilde{g}_{ab}^{-1}\circ \tilde{g}_{ab}: a^*V\to b^*V \to a^*V$ is the identity on $a^*(V^L)$ as $\tilde{g}_{ab}^{-1}\circ \tilde{g}_{ab}
\in L(T)_{aa}$. Thus the lifting $\tilde{g}_{ab}$ yields a well defined action of $G(K/k)$ on $H^0(L, V)$. 

One considers the following diagram of functors:
\ga{6.1}{\begin{CD}
\Rep(K: G(X/k))&@>{H^0(L,V)}>>& \Rep(K:G(K/k))\\
@VVV&&@VVV\\
\MIC(X/k)&@>>{H^0_{DR}(X, (\sV, \nabla/K))}>&\MIC(K/k)
\end{CD}}
According to Lemma \ref{lem:key}, the canonical
morphism
\ga{6.2}{H^0(L,V)\longrightarrow H^0_{DR}(X/K,(\sV,\nabla/K))}
is an isomorphism. Thus the above diagram is commutative. As a consequence
we obtain canonical morphisms
\ga{6.3}{\mathbf R^n_{G(X/k)}H^0(L,V)\longrightarrow
\mathbf R^n_{{\rm MIC}}H^0_{DR}(X,(\sV,\nabla/K))}
where, on the left hand-side, the derived functor is taken in
$\Rep(K: G(X/k))$ and
on the right hand-side the derived functor is taken in $\MIC(X/k)$.
From the Remark \ref{rmk4.1} and the commutative diagram \ref{4.3},
we know that the right hand-side is  the $n$-th relative de Rham cohomology
$$H^n_{DR}(X,(\sV, \nabla/K))=\mathbf R^n_{{\rm MIC}}H^0_{DR}(X,(\sV, \nabla/K))$$
equipped with the Gau\ss-Manin connection. It is a finite dimensional
$K$-vector space, as $\sV$ if of finite rank. Thus
$H^n_{DR}(X, (\sV, \nabla/K)) $ together with its Gau{\ss}-Manin connection
is an object of  $\Rep_f(K: G(K/k))$ and the homomorphism
in \eqref{6.3} is $G(K/k)$-equivariant.

The group cohomology $H^i(L,V)$ for any $L$-representation $V$
is defined as the right derived functor
of the functor $$\Rep(L)\xrightarrow{H^0(L, V)} {\rm Vec}_K^{qc}.$$
In case $V$ is the
restriction to $L$ of a representation of $G(X/k)$, there exists a
canonical homomorphism $$\mathbf R^n_{G(X/k)}H^0(L,V)\to H^n(L,V)$$
defined by constructing a map from an injective resolution of $V$
in $\Rep(K:G(X/k))$ to an injective resolution
 in $\Rep(L)$ (see Appendix, Lemma \ref{app:lem01}).

\begin{prop}\label{prop6.1} The canonical homomorphism
\ga{6.4}{\mathbf R^n_{G(X/k)}H^0(L,V)\to H^n(L,V)}
is an isomorphism. Consequently it induces a representation of $G(K/k)$
on $H^n(L, V)$ which has the property that the  canonical  homomorphism
\ga{6.5}{H^i(L,V)\to H^i_{ DR}(X,(\sV,\nabla/K))}
 is $G(K/k)$-equivariant.
\end{prop}
\begin{proof}
According to the discussion above, it suffices to show that
a representation of $G(X/k)$ which is injective (as an object in
$\Rep(K:G(X/k))$), remains injective when considered as a
representation of $L$.
The proof is based on the following lemma which will
be proved in the rest of the section.

Set $G:=G(X/k)$.
Let $\sO(G)$ be the ring of regular functions on $G$.
There is a natural action of $G$ on $\sO(G)$ called the
left regular action, see Appendix.
\begin{lem}\label{lem6.2} The $G$-representation $\sO(G)$ restricted to
$G^\Delta$ is injective in the category $\Rep(G^\Delta)$.
\end{lem}

Let us first assume this lemma. Since $L$ is normal in $G^\Delta$, $\sO(G^\Delta)$ is
injective as an $L$-representation. Indeed, $\sO(G^\Delta)$ is
faithfully flat over $\sO(G^\Delta/L)$, \cite{waterhause}, Chapter 16,
hence by \cite{takeuchi79}, Theorem 1, it is injective as an
$L$-representation.
Therefore any injective $G^\Delta$-representation, being direct summand
of a direct sum of copies of $\sO(G^\Delta)$, remains injective when considered as a
$L$-representation. Thus $\sO(G)$ is also an injective $L$-representation.

According to Lemma \ref{app:lem02} in the Appendix, we have
the following resolution of $V$ in $\Rep(K:G)$:
\ga{6.6}{\begin{CD}
V\otimes_t\sO(G)@>>> V
\otimes_t\sO(G)\otimes_t\sO(G)@>>> \ldots\\
@A\cong AA \\
V
\end{CD}}
where the tensor product is taken over $K$ and the action of $K$
on $\sO(G)$, indicated by the subscript ${}_t$, is induced from
the map $t:G\to \Spec(K)$. On each term
$V\otimes \sO(G)\otimes \sO(G)\ldots\otimes \sO(G)$ of the
complex, $G$ acts by its action
on the last tensor term. Hence,
as a complex of $G^\Delta$-modules, it is a resolution of
$V$ by injective $G^\Delta$-modules.

As \eqref{6.6} is an injective resolution of $V$ both in  $\Rep(K:G)$ and in
$\Rep(G^\Delta)$, its cohomology computes
$\mathbf R^*H^0(L,V)$ as well as $H^n(L,V)$.  This shows Proposition \ref{prop6.1}.
\end{proof}

The rest of this section is devoted to the proof of Lemma \ref{lem6.2}.

\subsection{The algebra $\sO(G^\Delta)$}
We refer to the Appendix for the properties of $\sO(G)$ and $\sO(G^\Delta)$.

By definition of $G^\Delta$, we have
\ga{6.7}{\sO(G^\Delta)\cong \sO(G)\otimes_{K\otimes_kK}K}
where $K\otimes_kK \to K$ is the product map.
 Then
$J:={\rm Ker}(K\otimes_kK\to K)$ is generated by elements
 of the form $\lambda\otimes 1
-1\otimes \lambda$, $\lambda\in K$.
 Since $\sO(G)$ is faithful over
$K\otimes_kK$, tensoring the exact sequence $0\to J\to  K\otimes_kK\to K \to 0$
with $\sO(G)$ over $ K\otimes_kK$, one obtains an exact sequence
\ga{6.8}{0\to J\sO(G)\to \sO(G)\xrightarrow{\pi} \sO(G^\Delta) \to 0.}
That is, we can identify $J\otimes_{K\otimes_kK}\sO(G)$ with
its image $J\sO(G)$ in $\sO(G)$.

\subsection{The functor $\Ind$}
For any representation $W\in\Rep(G^\Delta)$,
define $\Ind(W)$ to be
\ga{6.9}{\Ind(W):=(W\otimes_t\sO(G))^{G^\Delta}}
where $G^\Delta$ acts on $W$ as usual and on $\sO(G)$ through
the right regular action of $G$ on $\sO(G)$ (i.e $\sO(G)$
is a right $G$-module). On this invariant space, $G$ acts through
the left regular action on $\sO(G)$. Thus $\Ind$ is a functor
$\Rep(G^\Delta)\to \Rep(K:G)$.

The space $\Ind(W)$ can also be given as the equalizer of the maps
\ga{6.10}{\begin{array}{l} p: W\otimes_t\sO(G)
\stackrel{\rho_W\otimes\id}\longrightarrow
W\otimes\sO(G^\Delta)\otimes_t\sO(G)\\
q:W\otimes_t\sO(G)\stackrel{\id\otimes\Delta}\longrightarrow
 W\otimes_t\sO(G)\stotimes\sO(G)\stackrel\pi\longrightarrow
W\otimes\sO(G^\Delta)\otimes_t\sO(G)\end{array}}
where $\rho_W:W\to W\otimes \sO(G^\Delta)$ is the coaction of
$\sO(G^\Delta)$ on $W$, $\Delta$ is the coproduct on $\sO(G)$.

\begin{lem}\label{lem6.3}
There exists a functorial isomorphism
\ga{6.11}{\Hom_G(V,\Ind (W))\cong \Hom_{G^\Delta}(V,W),}
$V\in \Rep(K:G), W\in\Rep(G^\Delta)$,
i.e., $\Ind$ is the right adjoint to the functor restricting
$G$-representations to $G^\Delta$.
\end{lem}
\begin{proof} The map is given by composing with the canonical
projection $\Ind (W)\to W$, $v\otimes h\mapsto v\varepsilon(h)$. The
converse map is given by $f\mapsto (f\otimes\id)\rho_W$.
\end{proof}

\begin{lem}\label{lem6.4} The functor $\Ind$ is exact if and only if
$\sO(G)$
is injective as a $G^\Delta$-module.
\end{lem}
\begin{proof}
Since $G^\Delta$ is a group scheme over a field $K$,
its representations are union of their subrepresentations
 of finite dimension over
$K$. Therefore the injectivity of $\sO(G)$ requires only to be
checked on finite dimensional representations of $G^\Delta$.
For such a representation $W$, we have
$$\Hom_{G^\Delta}(W,\sO(G))\cong\Hom_{G^\Delta}(K,W^*\otimes\sO(G))=
\Ind(W^*)$$
Since the dualizing functor $(-)^*$ and the functor tensoring over $K$
are exact, the claim follows.
\end{proof}

Let us use the following notation of Sweedler for the coproduct on
$\sO(G)$:
$$\Delta(g)=\sum_{(g)}g_{(1)}\otimes g_{(2)}.$$

\begin{lem}\label{lem6.5} The following map
\ga{6.12}{\varphi:\sO(G)\otimes_{K\otimes_kK}\sO(G)\to
\sO(G^\Delta)\otimes_K{}_t\sO(G)\\
g\otimes h\mapsto \sum_{(g)}\pi(g_{(1)})\otimes g_{(2)}h \notag}
is an isomorphism, where $\pi$ is defined  in formula \eqref{6.8}.
\end{lem}
\begin{proof} We define the inverse map to this map. Let
\ga{6.13}{\bar\psi:\sO(G)\stotimes \sO(G)\to \sO(G)\otimes_{K\otimes_kK}
\sO(G)}
be the map that maps $g\otimes h\mapsto \sum_{(g)}g_{(1)}\otimes
\iota(g_{(2)}h$. We have for $\lambda\in K$, and for $t,s:K\to\sO(G)$
\begin{eqnarray*}
\bar\psi(t(\lambda)g\stotimes h)&=&
\sum_{(g)} g_{(1)}\otimes \iota(t(\lambda)g_{(2)})h\\
&=&\sum_{(g)} g_{(1)}\otimes_{K\otimes_kK} s(\lambda)\iota(g_{(2)})h
 \quad\mbox{ by (\ref{A13}) }\\
&=&s(\lambda)\sum_{(g)} g_{(1)}\otimes_{K\otimes_kK} \iota(g_{(2)})h\\
&=&\bar\psi(s(\lambda)g\stotimes h)
\end{eqnarray*}

Thus $\bar\psi$ maps $J\sO(G)\stotimes \sO(G)$ to 0,
hence factors through a map $\psi:\sO(G^\Delta)\otimes_t\sO(G)\to
\sO(G)\otimes_{K\otimes_kK}\sO(G)$. Checking $\varphi\psi=\id$,
$\psi\varphi=\id$ can be easily done using the property
(\ref{A14}) of $\iota$.
\end{proof}

\begin{cor}\label{cor6.6} For any $W\in \Rep(G^\Delta)$ we have the following isomorphism
\ga{6.14}{\Phi:\Ind (W)\otimes_{K\otimes_kK}\sO(G)\cong W\otimes_t\sO(G)\\
\Phi(w\otimes g\otimes h)=w\otimes gh.\notag} The inverse is given
by 
\ga{}{
W\otimes_t\sO(G)\to W\otimes \sO(G^\Delta)\otimes\sO(G)
\xrightarrow{\Psi} W\otimes_t\sO(G)\otimes_{K\otimes_kK}\sO(G)\notag\\
\Psi= (\id_W\otimes\psi)(\rho_W\otimes\id).\notag
}
\end{cor}
\begin{proof} Tensoring the isomorphism in (\ref{6.12}) with
$W$ and applying the functor $(-)^{G^\Delta}$ we obtain
$\Phi$.
\end{proof}

\subsection{Proof of Lemma \ref{lem6.2}} According to Lemma \ref{lem6.4}, it
suffices to show the exactness of $\Ind$. According to Corollary \ref{cor6.6},
the functor
$$\Ind(-)\otimes_{K\otimes_kK}\sO(G)\cong (-)\otimes_t\sO(G)$$
hence is exact. Since $\sO(G)$ is faithfully flat over ${K\otimes_kK}$,
$\Ind$ is exact. \hfill $\Box$
\begin{rmk}\label{rmk6.7} The above proof for groupoid schemes
is inspired by  Takeuchi's proof (\cite{takeuchi79}) for
the case of group schemes.
\end{rmk}

\begin{appendix}\section{Groupoid schemes}

In this Appendix, we briefly recall  the notions of affine groupoids and
their representations which are used in the article. Our reference
is \cite{DeGroth}, Section 3.
\subsection{Groupoid schemes} We fix a field $k$. By a $k$-affine
scheme we mean the spectrum of a $k$-algebra (not necessarily
finitely generated over $k$). 
Let $S/k$ be a $k$-affine scheme.
 With this terminology, $S$ can be taken to be the spectrum of a field extension $S=K\supset k$. 
An affine $k$-groupoid scheme acting on $S$
is a $k$-affine scheme $G$ together with
two morphisms (the source and the target maps) $s,t:G\to S$, satysfying the
following axioms:
\begin{itemize}
\item[(i)] There exists a map $m:G\sttimes G\to G$ called the product of
$G$, satisfying the following associativity property
\ga{A1}{m(m\sttimes\id_G)=m(\id_G\sttimes m)}
\item[(ii)] There exists a map $\varepsilon:S\to G$ called the unit
element map, satisfying the following property
\ga{A2}{m(\varepsilon\sttimes\id_G)=m(\id_G\sttimes m)=\id_G}
\item[(iii)] There exists a map $\iota:G\to G$, called the inverse map,
satisfying the following properties:
\ga{A3}{\iota\circ s=t;\quad \iota\circ t=s}
\ga{A4}{m(\iota \sttimes \id_G)=\varepsilon\circ s,\quad
m(\id_G\sttimes\iota)=\varepsilon\circ t,}
\end{itemize}
where $\sttimes$ denotes the fiber product over $S$ with respect
to the maps $s$ and $t$.

Let $T$ be a $k$-schema. By definition, the category $(S(T),G(T))$ has for 
objects the morphisms $T\to S$ and for morphisms between
two objects $a,b:T\to S$ the morphisms $\phi:T\to G$
satifying
\ga{A5}{(a,b)=(s,t)\phi:T\to S\times S.}
The axioms A.1--4 for $G$ imply that this category is
a groupoid.
Note that the set of all  morphisms of
$(S(T),G(T))$ is precisely $G(T)=\Hom_k(T,G)$.

A groupoid scheme $G$ acting on $S$ is said to be acting transitively
if there is a  map $\phi:T\to G$ such that the
map $(s,t)\phi:T\to S\times S$ is faithfully flat.
This implies that $(s,t):G\to S\times S$
is a faithfully flat map.

A groupoid scheme $G$ acting on $S$ is called discrete if
the structure map $(s,t)$ factors through the diagonal
map $\Delta:S\to S\times S$  and a map $u:G\to S$.
In this case $G$ equipped with  $u$ is an $S$-group scheme.

For example, define 
$G^\Delta$ as the
pull-back  of $G$ along the diagonal map $\Delta:S\times S$.
\ga{A6}{
\begin{CD}
G^\Delta @>>> G\\
@V(s,t)VV @VV(s,t)V\\
\Delta @> {\rm diagonal}>>  S\times_k S \end{CD}
}
Then $G^\Delta$ is a discrete $S$-groupoid scheme, which is
a subgroupoid scheme in $G$.

Another simple example is $S$, which is a groupoid acting on
itself by means of the diagonal map.
\subsection{Homomorphisms}
A morphims of $k$-groupoid schemes acting on a $k$-scheme
$S$ is a morphism of the underlying $k$-schemes which is
compatible with all structure maps. For instance, the
unit element map $\varepsilon: S\to G$ is a morphism
of groupoid schemes.

For two homomorphisms of groupoid schemes $G_i\to G$,
 $i=1,2$, there exists an obvious structure of
 groupoid scheme on $G_1\times_G G_2$. In particular,
 we define the kernel of a homomorphism $f:G_1\to G$
 as the fiber product $\ker f:=S\times_GG_1$. It is easy to see
 that $\ker f$ is a discrete groupoid scheme, defining a group scheme
 over $S$. Assuming that $G_1$ and $G$ act transively on $S$, then,
 by taking the fiber product with $S$ over $S\times S$,
 i.e. taking the diagonal group schemes, we see
 that $\ker f$ is isomorphic to the kernel of the  homomorphism
 $G_1^\Delta\to G^\Delta$ of group schemes:
 \ga{A7}{
 \begin{CD} \ker f&@>>>&G_1&@>{f}>>& G\\
 @|& & @AAA& & @AAA\\
  \ker f^\Delta&@>>>& G_1^\Delta&@>>{f^\Delta}>&G\end{CD}}

\subsection{Representation} Let $V$ be a quasi-coherent
sheaf on $S$. A representation of $G$ in $V$ is an operation
 $\rho$, that assigns to each $k$-schema $T$ and each
 morphism $\phi:T\to G$ a $T$-isomorphism
\ga{A8}{\rho(\phi):a^*V\to b^*V}
where $(a,b)=(s,t)\phi$, the source and the target of $\phi$, and
$a^*$ (resp. $b^*$) denotes the pull-back of $V$
along $a$ (resp. $b$). One requires that this operation
be compatible with the composition law of the groupoid
$(S(T),G(T))$ and with the base change.  The latter means:
for any morphism $r:T'\to T$
\ga{A9}{\rho(r^*\phi)=r^*\rho(\phi).}

In particular, one has the trivial representation of $G$ in
$R=\sO_S$ where all morphisms $\rho(\phi)$ are identity morphisms.

\subsection{Tannaka duality}\label{sub:tannaka}
 Assume that $G$ acts transitively on $S$, then
 representations of $G$ form
an abelian category which is closed under taking the
tensor product. We denote this category by
$\Rep(S:G)$. We denote the full subcategory of $\Rep(S:G)$
of those representations which are of finite rank as
sheaf on $S$ by $\Rep_f(S:G)$. Each object of $\Rep_f(S:G)$
is local free when considered as sheaf on $S$,  and each object
of $\Rep(S:G)$ is a filtered union of its finite rank subrepresentations.
Using the inverse map, to each representation
in a coherent locally free $\sO_S$-module, 
one can define a representation in the dual coherent sheaf.
Finally, for the trivial representation
in $\sO_S$, the set of endomorphisms is isomorphic to $k$.
See \cite{DeGroth}, Section 3, for details.

A category with the above properties is called
a tensor category over $k$. Conversely, for any
tensor category $\sC$ over $k$ with a fiber functor
to the category Qcoh$(S)$ of quasi-coherent sheaves over $S$, one
can construct a groupoid scheme $G$ acting transitively on $S$,
such that the fiber functor factors becomes an equivalence
of tensor categories $\sC\stackrel\cong\to \Rep(S:G)$. This correspondence
is in fact a 1-1 correspondence between tensor categories
over $k$ equipped with a fiber functor to Qcoh($S$) and
$k$-groupoids acting transitively over $S$, known as
the Tannaka duality \cite{DeGroth}, Th\'eor\`eme 1.12.

\subsection{Representations of discret groupoids}
If $G$ is a groupoid scheme acting discretly over $S$ then one can
easily deduce from the definion that representations of $G$
are in 1-1 correspondence with representations of the
underlying $S$-group scheme.
If $(\rho,V)$ is a  representation, for all
commutative diagrams
$$\begin{CD}
T @>g_{(a,a)}>> G\\
@V=VV @VVV \\
T @> (a,a)>> \Delta=S
         \end{CD}$$
 with $T$ a $k$-scheme, then $\rho(g_{a,a})$ is an isomorphism from $a^*V$ to itself. But  $ G(T)=\emptyset$ for $a\neq b$.
In other words, $\rho$ induces a representation of $G$, as an $S$-group scheme, in $V$. The converse is also true.

\subsection{The function algebra} Let $R:=\sO(S)$ denote the algebra of
regular functions on $S$.
The groupoid structure on $G$ induces the following structures on $\sO(G)$.
The source and the target map for $G$ induce algebra maps
$s,t:K\to \sO(G)$.
The transitivity of $G$ on $S$ can be rephrased by saying that
$\sO(G)$ is faithfully flat over $R\otimes_kR$ with
respect to the base map $t\otimes_ks:R\otimes_kR\to \sO(G)$.

The composition law for $G$ induces an $R\otimes_kR$-algebra map
\ga{A10}{\Delta: \sO(G)\longrightarrow \sO(G)\stotimes  \sO(G).}
satisfying $(\Delta\otimes \id)\Delta=(\id\otimes \Delta)\Delta$.
The unit element of $G$ induces a $R\otimes_kR$-algebra map
\ga{A11}{\varepsilon:\sO(G)\longrightarrow R}
where $R\otimes_kR$ acts on $R$ diagonally (i.e., $\lambda\otimes_k\mu)\nu
=\lambda\mu\nu$).
One has
\ga{A12}{(\varepsilon\otimes \id)\Delta=(\id\otimes\varepsilon)\Delta
=\id}

Finally, the operation which consists of taking the inverse  in $G$ induces an automophism
$\iota$ of $\sO(G)$ which interchanges the actions $t$ and $s$:
\ga{A13}{\iota(t(\lambda)s(\mu)h)=s(\lambda)t(\mu)\iota(h)}
and satisfies the following equations:
\ga{A14}{ m(\iota\otimes\id)\Delta=s\circ\varepsilon\quad
m(\id\otimes\iota)\Delta=t\circ\varepsilon
}

A representation $\rho$ of $G$ in $V$ induces a map
$\rho:V\to V\otimes_t \sO(G)$, called coaction of $\sO(G)$ on $V$, such that
\ga{A15}{(\id_V\otimes \Delta)\rho=(\rho\otimes \id_V),\quad
(\id_V\otimes\varepsilon)\rho=\id_V.}
An $R$-module equipped with such an action is called $\sO(G)$-comodule.
Conversely, any coaction of $\sO(G)$ on an $R$-module $V$ defines
a representation of $G$ in $V$. In fact, we have an equivalence
between the category of $G$-representations and the category
of $\sO(G)$-comodules. The discussion in the previous
subsection shows that $V$ is projective over $R$.

In particular, the coproduct on $\sO(G)$ can be considered as a coaction
of $\sO(G)$ on itself and hence defines a representation of $G$ in
$H$, called the right regular representation.
\begin{lem}\label{app:lem01}
Let $G$ be a $k$-groupoid scheme acting transitively on $S$.
Consider the function algebra $\sO(G)$ as $G$-representation
with respect to the left regular action.
Then for any injective $R$-module $V$, the $G$-representation
$V\otimes_R\sO(G)$, where the action of $G$ is given by the action of
 $\sO(G)$, is an injective object in $\Rep(R:G)$. In particular, when $R=K$ is a field, 
${\rm Rep}(K: G)$ has enough injective objects. 
\end{lem}
\begin{proof} For any $\sO(G)$-comodules $U$
we have the following functorial isomorphism
\ga{A16}{\Hom_G(U,V\otimes_t \sO(G))\cong \Hom_R(U,V)\\
f\longmapsto (\id\otimes\varepsilon)f.\notag} Indeed, the inverse is given
by
$$g\longmapsto \rho_V\circ g$$ where
$\rho_V$ is the coaction of $\sO(G)$ on $V$.

Finally, if $K$ is a field, all $K$-vector spaces are injective objects. 
\end{proof}

\begin{lem}\label{app:lem02}
Assume that $R=K$ is a field. Let $V\in \Rep(K:G)$.
Then the following complex is a resolution
of $V$ in $\Rep(S:G)$.
\ga{A17}{\begin{CD}
 V\otimes_t \sO(G)@>>> V
\otimes_t\sO(G)\otimes_t\sO(G)@>>> \ldots \\
@A\cong AA \\
V
\end{CD}}
where the tensor product is taken over $K$ and the index $t$
specifies the action $t$ of $K$.
\end{lem}
\begin{proof} This can be done exactly as in the case of
group schemes over a field (\cite{Jan}, Chapter 4),
and will be omitted.\end{proof}

\end{appendix}

\bibliographystyle{plain}

\renewcommand\refname{References}

\end{document}